\documentclass [12pt,a4paper,reqno]{amsart}
\usepackage{amsmath,amsthm}
\usepackage{amsfonts}
\usepackage{amssymb}
\usepackage{tikz-cd}
\allowdisplaybreaks

\DeclareMathOperator{\esssup}{ess\,sup}

\newcommand{\be}{\begin{equation}}
\newcommand{\ef}{\end{equation}}
\chardef\bslash=`\\ % p. 424, TeXbook
\newtheorem*{thm*}{Theorem}

\theoremstyle{definition}

\newtheorem*{remark*}{Remarks}
\newtheorem*{defn*}{Definition}

\theoremstyle{remark}
\newcommand{\wh}{\widehat}
\newcommand{\fc}{\frac}

\newcommand{\iy}{\infty}
 \renewcommand{\sectionmark}[1]{}
\renewcommand{\Im}{\operatorname{Im}}

\newcommand{\ve}{\varepsilon}

\newcommand{\const}{\operatorname{const}}

 \usepackage{amsfonts}
\newcommand{\field}[1]{\mathbb{#1}}

\newcommand{\D}{\field{D}}
\newcommand{\om}{\omega}
\newcommand{\z}{\zeta}
\newcommand{\ov}{\overline}
\newcommand{\vp}{\varphi}
\newcommand{\hC}{\wh{\field{C}}}
\newcommand{\C}{\field{C}}
\newcommand{\R}{\field{R}}

\newcommand{\B}{\mathbf{B}}
\newcommand{\T}{\mathbf{T}}
\newcommand{\Belt}{\operatorname{Belt}}

\newcommand{\vk} {\varkappa}
\newcommand{\x} {\mathbf x}
\renewcommand{\a} {\alpha}

\newcommand{\ld}{\lambda}

\begin{document}

\title{Extremal properties of Sobolev's Beltrami coefficients and distortion of curvelinear functionals}
\author{Samuel L. Krushkal}

\begin{abstract}
 An important problem in applications of quasiconformal analysis and in its numerical aspect is to establish algorithms
for explicit or approximate determination of the basic quasiinvariant curvelinear and analytic functionals intrinsically
connected with conformal and quasiconformal maps, such as their Teichm\"{u}ller and Grunsky norms, Fredholm eigenvalues and the
quasireflection coefficients of associated quasicircles.

We prove a general theorem of new type answering this question for univalent functions in arbitrary quasiconformal domains and provide its applications. The results are strengthened in the case of maps of the disk and give rise to extremal Beltrami coefficients of a new type.
\end{abstract}

\date{\today\hskip4mm({SobolBel1.tex})}

\maketitle

\bigskip

{\small {\textbf {2020 Mathematics Subject Classification:} Primary: 30C62, 30C75; Secondary: 30F60, 32F45,32G15, 46G20}

\medskip

\textbf{Key words and phrases:} Univalent function, Grunsky operator, quasiconformal extension, the Sobolev spaces,
harmonic Beltrami coefficients, quasicircles, quasireflections, universal Teichm\"{u}ller space}

\bigskip

\markboth{Samuel L. Krushkal}{Extremal properties of Sobolev's Beltrami coefficients} \pagestyle{headings}

\bigskip\bigskip
\centerline{\bf 1. PRELIMINARIES}

\bigskip\noindent
{\bf 1.1. Preamble}.
 An important still open problem in geometric complex analysis is to establish algorithms for explicit
or approximate determination of the basic curvilinear and analytic functionals intrinsically connected with
conformal and quasiconformal maps, such as their Teichm\"{u}ller and Grunsky norms, Fredholm eigenvalues and the
quasireflection coefficients of associated quasicircles. It is important also for the potential theory.
However, this has not been solved completely even for convex polygons.

The problem has intrinsic interest also in view of its connection
with geometry of Teichm\"{u}ller spaces and with the approximation theory. It is crucial also for numerical aspect
of quasiconformal analysis.

This paper deals with functions not admitting the canonical Teichm\"{u}ller extensions and solves explicitly the indicated
problem for a natural rather broad class of Beltrami coefficients supported in generic quasiconformal domains,
in other words, for univalent functions with quasiconformal extension on arbitrary quasidisks.
The results give rise to extremal Beltrami coefficients of a new type.

\bigskip\noindent
{\bf 1.2. The Teichm\"{u}ller and Grunsky norms of univalent functions}.
We start with the class $\Sigma_Q$ of univalent functions $f(z) = b_0 + b_1 z^{-1} + \dots$ in the disk
$\D^* = \{z \in \hC: \ |z| > 1\}$ admitting quasiconformal extensions across the boundary unit circle $\mathbb S^1
= \partial \D^*$, hence to the whole Riemann sphere $\hC = \C \cup \{\iy\}$. To have compactness in the topology of
locally uniform convergence on $\C$, one must add the third normalization condition, for example,
$f(0) = 0$.

The Beltrami coefficients of extensions are supported in the unit disk
$\D = \{|z| < 1\}$ and run over the unit ball
$$
\Belt(\D)_1 = \{\mu \in L_\iy(\C): \ \mu(z)|\D^* = 0, \ \ \|\mu\|_\iy  < 1\}.
$$
Each $\mu \in \Belt(\D)_1$ determines a unique homeomorphic solution to the Beltrami equation
$\ov \partial w = \mu \partial w$ on $\C$ (quasiconformal automorphism of $\hC$)  normalized by
the assumptions $w^\mu \in \Sigma_Q, w^\mu(0) = 0$.

The {\bf Schwarzian derivatives} of these functions
$$
S_w(z) = \Bigl(\fc{w^{\prime\prime}(z)}{w^\prime(z)}\Bigr)^\prime
- \fc{1}{2} \Bigl(\fc{w^{\prime\prime}(z)}{w^\prime(z)}\Bigr)^2, \quad z \in \D^*
$$
belong to the complex Banach space $\B = \B(\D^*)$ of hyperbolically bounded holomorphic functions in
the disk $\D^*$ with norm
$$
\|\vp\|_\B = \sup_\D \ (|z|^2 - 1)^2 |\vp(z)|
$$
and run over a  bounded domain in $\B$ modeling the {\bf universal Teichm\"{u}ller space} $\T$.
The space $\B$ is dual to the Bergman space $A_1(\D^*)$, a subspace of $L_1(\D^*)$ formed by integrable holomorphic functions (quadratic differentials $\vp(z) dz^2$) on $\D^*$. Note that $\vp(z) = O(z^{-4})$ near $z = \iy$.
The needed results from Teichm\"{u}ller space theory see, e.g., in \cite{EKK}, \cite{GL}, \cite{Kr1}.

One defines for any $f \in \Sigma_Q$ its {\bf Grunsky coefficients} $\a_{m n}$ from the expansion
 \be\label{1}
\log \fc{f(z) - f(\z)}{z - \z} =  \sum\limits_{m, n = 1}^\iy \a_{m n} z^{-m} \z^{-n}, \quad (z, \z) \in (\D^*)^2,
\end{equation}
where the principal branch of the logarithmic function is chosen. These coefficients satisfy
the inequality
 \be\label{2}
\Big\vert \sum\limits_{m,n=1}^\iy \ \sqrt{m n} \ \a_{m n} x_m x_n \Big\vert \le 1
\end{equation}
for any sequence $\mathbf x = (x_n)$ from  the unit sphere $S(l^2)$ of the Hilbert space $l^2$
with norm $\|\x\| = \bigl(\sum\limits_1^\iy |x_n|^2\bigr)^{1/2}$; conversely, the inequality (2)
also is sufficient for univalence of a locally univalent function in $\D^*$ (cf. \cite{Gr}, \cite{Mi}).

The minimum $k(f)$ of dilatations $k(w^\mu) = \|\mu\|_\iy$ among all quasiconformal extensions
$w^\mu(z)$ of $f$ onto the whole plane $\hC$ (forming the equivalence class of $f$) is called the \textbf{Teichm\"{u}ller norm} of this function. Hence,
$$
k(f) = \tanh d_\T(\mathbf 0, S_f),
$$
where $d_\T$ denotes the Teichm\"{u}ller-Kobayashi metric on $\T$. This quantity dominates the \textbf{Grunsky norm}
$$
\vk(f) = \sup \Big\{\Big\vert \sum\limits_{m,n=1}^\iy \ \sqrt{mn} \
\a_{m n}(f) x_m x_n \Big\vert: \ \mathbf x = (x_n) \in S(l^2) \Big\}
$$
by $\vk(f) \le k(f)$.
For most functions $f$, we have the strong inequality $\vk(f) < k(f)$ (moreover, the functions satisfying this
inequality form a dense subset of $\Sigma$), while the functions with the equal norms play a crucial role in many
applications.

These norms coincide only when any extremal Beltrami coefficient $\mu_0$ for $f$ (i.e., with $\|\mu_0\|_\iy = k(f)$)
satisfies
 \be\label{3}
\|\mu_0\|_\iy = \sup \ \Big\{ \Big\vert \iint_\D \mu_0(z) \psi(z) dx dy \Big\vert: \ \psi \in A_1^2(\D), \ \|\psi\|_{A_1(\D)} = 1\Big\} = \vk(f) \quad (z = x + i y).
\end{equation}
Here $A_1(\D)$ denotes the subspace in $L_1(\D)$ formed by integrable holomorphic functions
(quadratic differentials $\psi(z) dz^2$ on $\D$, and $A_1^2(\D)$ is its subset consisting of $\psi$ with zeros of even
order on $\D$, i.e., of the squares of holomorphic functions (see, e.g., \cite{Kr2}, \cite{Kr4}, \cite{Kr7}, \cite{Ku2}).
Note that every $\psi \in A_1^2(\D)$ has the form
 \be\label{4}
\psi(z) = \fc{1}{\pi} \sum\limits_{m+n = 4}^\iy \sqrt{m n} \ x_m x_n z^{m+n-2}
\end{equation}
and $\|\psi\|_{A_1(\D)} = \|\x\|_{l^2} = 1, \ \x = (x_n)$.

\bigskip\noindent
{\bf 1.3. Generalization}.
The method of Grunsky inequalities was generalized in several directions, even to bordered Riemann surfaces $X$ with a finite number of boundary components (see, e.g., \cite{Mi}, \cite{SS}).
We shall consider these inequalities in unbounded simply connected hyperbolic domains.

Let $L \subset \C$ be an oriented bounded quasicircle separating the points $0$ and $\iy$. Denote its interior and exterior domains by $D$ and $D^*$ (so $0 \in D, \ \iy \in D^*$).
Then, if $\delta_D(z)$ denotes the Euclidean distance of $z$ from the boundary of $D$ and $\ld_D(z) |dz|$ is its hyperbolic
metric of Gaussian curvature $-4$, we have
 \be\label{5}
\fc{1}{4} \le \ld_D(z) \delta_D(z) \le 1.
\end{equation}
The right hand inequality follows from the Schwarz lemma and the left from the Koebe one-quoter theorem.

For such a domain $D^* \ni \iy$, the expansion (1) assumes the form
$$
- \log \fc{f(z) - f(\z)}{z - \z} = \sum\limits_{m, n = 1}^\iy
\fc{\beta_{m n}}{\sqrt{m n} \ \chi(z)^m \ \chi(\z)^n},
$$
where $\chi$ denotes a conformal map of $D^*$ onto the disk $\D^*$ so that $\chi(\iy) = \iy, \ \chi^\prime(\iy) > 0$ (cf. \cite{Mi}).

Accordingly, the generalized Grunsky norm is defined by
$$
\vk_{D^*}(f) = \sup \Big\{ \Big\vert \sum\limits_{m,n = 1}^{\iy} \ \beta_{mn} \ x_m x_n \Big\vert : \
{\mathbf x} = (x_n) \in S(l^2)\Big\}.
$$

We now consider the class $\Sigma_Q(D^*)$ of univalent functions in domain $D^*$ with expansions $f(z) = b_0 + b_1 z^{-1}
+ \dots$ near $z = \iy$, admitting quasiconformal extensions onto the complementary domain $D$. Similar to above, we subject these extensions to $f(0) = 0$.
Their Beltrami coefficients run over the ball
$$
\Belt(D)_1 = \{\mu \in L_\iy(\C): \ \mu(z)|D^* = 0, \ \ \|\mu\|_\iy  < 1\}.
$$
A coefficient $\mu_0 \in \Belt(D^*)_1$ is {\it extremal in its class if and only if}
$$
\|\mu_0\|_\iy = \sup \ \Big\{ \Big\vert \iint_{D^*} \mu_0(z) \psi(z) dx dy \Big\vert: \ \psi \in A_1(D^*), \ \|\psi\|_{A_1} = 1\Big\},
$$
and similar to (3), the equality $\vk_{\D^*}(f^\mu) = k(f^\mu)$ is valid if and only if
$$
\|\mu\|_\iy = \sup \ \Big\{ \Big\vert \iint_D \mu(z) \psi(z) dx dy \Big\vert: \ \psi \in A_1^2(D), \ \|\psi\|_{A_1(D)} = 1\Big\} = \vk_D(f).
$$
If additionally the equivalence class of $f$ is a {\bf Strebel point} of the space $\T$ with base point $D^*$,
which means that this class contains the Teichm\"{u}ller extremal extension  $f^{k|\psi_0|/\psi_0}$
with $\psi_0 \in A_1(D)$, then necessarily
 $\psi_0 = \om ^2 \in A_1^2$ (cf. \cite{Kr2}, \cite{Kr7}, \cite{Ku3}), \cite{St2}).
The Strebel points are dense in any Teichm\"{u}ller space, see \cite{GL}.

For arbitrary quasidisks $D$, the corresponding set $A_1^2(D)$ is characterized similar to (4), but in more complicated way
(see \cite{Kr8}).

\bigskip
Assume that $\mu_0 \in \Belt(D)_1$ is extremal in its class but not of Teichm\"{u}ller type.
A point $z_0 \in \partial D$ is called {\bf substantial} (or essential) for $\mu_0$ if
for any $\ve > 0$ there exists a neighborhood $U_0$ of $z_0$ such that
$$
\sup_{D^*\setminus U_0} |\mu_0(z)| < \|\mu_0\|_\iy - \ve;
$$
so the maximal dilatation $k(w^{\mu_0}) = \|\mu\|_\iy$ is attained on $D$ by approaching
this point.

In addition, there exists a sequence $\{\psi_n\} \subset A_1(D)$ such that $\psi_n(z) \to 0$
locally uniformly on $D$ but $\|\psi_n\| = 1$ for any $n$, and
$$
\lim\limits_{n\to \iy} \iint\limits_D \mu_0(z) \psi_n(z) dx dy = \|\mu_0\|_\iy.
$$
Such sequences are called {\bf degenerated}.

The image of a substantial point is a common point of two quasiconformal arcs, which can be of
spiral type.

\bigskip\noindent
{\bf 1.4. Fredholm eigenvalues and quasireflections}.
The Teichm\"{u}ller and Grunsky norms are intrinsically connected with quasiconformal reflections, Fredholm eigenvalues and other quasiinvariants of quasiconformal curves. We outline briefly the main notions; the details  see, e.g., in \cite{Ah2},  \cite{Kr3}, \cite{Kr8}, \cite{Ku4}.

\bigskip
The {\bf quasiconformal reflections} (or quasireflections) are the orientation reversing quasiconformal homeomorphisms of the sphere $\hC$ which preserve point-wise some (oriented) quasicircle $L \subset \hC$ and interchange its interior and exterior domains.
One defines for $L$ its {\bf reflection coefficient}
$$
 q_L = \inf k(f) = \inf \ \| \partial_z f/\partial_{\ov z} f \|_\iy,
$$
taking the infimum over all quasireflections across $L$. Due to \cite{Ah2}, \cite{Ku4}, the dilatation
$$
Q_L = (1 + q_L)/(1 - q_L) \ge 1
$$
is equal to the quantity
 \be\label{6}
Q_L = (1 + k_L)^2/(1 - k_L)^2,
\end{equation}
where $k_L$ is the minimal dilatation among all orientation preserving quasiconformal automorphisms $f_{*}$ of $\hC$ carrying the unit circle onto $L$, and
$k(f_{*}) = \|\partial_{\ov z} f_{*}/\partial_z f_{*}\|_\iy$.

The reflection with dilatation $Q_L$ is extremal. A remarkable and very useful fact established by Ahlfors is that any quasicircle also admits a Lipschitz continuous quasireflection  with some coefficient $C(q_L)$ (see \cite{Ah2}).

\bigskip
The {\bf Fredholm eigenvalues} $\rho_n$ of an oriented smooth closed Jordan curve $L \subset \hC$ are the eigenvalues of its double-layer potential. These values are crucial in many applications.

The least positive eigenvalue $\rho_L = \rho_1$ plays is naturally connected with conformal and quasiconformal maps and can be defined for any oriented closed Jordan curve $L$ by
$$
\fc{1}{\rho_L} = \sup \ \fc{|\mathcal D_G (u) - \mathcal D_{G^*}
(u)|} {\mathcal D_G (u) + \mathcal D_{G^*} (u)},
$$
where $G$ and $G^*$ are, respectively, the interior and exterior of
$L; \ \mathcal D$ denotes the Dirichlet integral, and the supremum
is taken over all functions $u$ continuous on $\hC$ and harmonic on
$G \cup G^*$.

A rough upper bound for $\rho_L$ is given by Ahlfors' inequality
$$
\fc{1}{\rho_L} \le q_L,
$$
where $q_L$ denotes the minimal dilatation of quasireflections across $L$ \cite{Ah1}.

One of the basic tools in quantitative estimating the Freholm eigenvalues $\rho_L$ of quasicircles is given by the
K\"{u}hnau-Schiffer theorem \cite{Ku2}, \cite{Sc}, which states that {\it the value $\rho_L$ is reciprocal
to the Grunsky norm $\vk(f)$ of the Riemann mapping function of the exterior domain of $L$}.

For all functions $f \in S_Q$ (i.e., univalent in the disk $\D^*$) with $k(f) = \vk(f)$, we have the exact explicit values
 \be\label{7}
q_{f(\mathbb S^1)} = \fc{1}{\rho_{f(\mathbb S^1)}} = \vk(f).
\end{equation}

\bigskip\noindent
{\bf 1.5. Harmonic and pseudo-harmonic Beltrami coefficients}.
By the Ahlfors-Weill theorem strengthening the classical Nehari's result on univalence in
terms of the Schwarzians, every function $\vp \in \B(\D)$ with $\|\vp\|_\B < 2$ is the Schwarzian derivative
of a univalent function $f(z)$ in the unit disk $\D$, and the function $f$ has
quasiconformal extension onto the disk $\D^*$ with Beltrami coefficient
  \be\label{8}
\mu_\vp(z) = - \frac{1}{2}(|z|^2 - 1)^2 \vp(1/{\ov z})(1/{\ov z}^4), \quad z \in \D^*;
\end{equation}
see \cite{AW}, \cite{Be}, \cite{GL}.

This deep fact is extended to arbitrary quasidisks in much weaker form. For example, we have the following result of Bers,
which is a special case of his more general extension theorem given in \cite{Ber}.

\bigskip\noindent
{\bf Lemma 1}. \cite{Ber} {\it Let $L$ be a quasicircle on $\hC$ with the interior $D_L$ and exterior $D_L^*$.
Then, for some $\ve > 0$, there exists an anti-holomorphic homeomorphism $\tau$ (with $\tau( \mathbf 0) = \mathbf 0$) of the ball \ $V_\ve = \{\vp \in \B(D_L^*): \|\vp\|\} < \ve$ \ into $\B(D_L)$ such that every $\vp$ in $V_\ve$ is the Schwarzian derivative of some univalent function $f$ which is the restriction to $D_L^*$ of a quasiconformal automorphism $\wh f$ of Riemann sphere $\hC$. This $\wh f$ can be chosen in such a way that its Beltrami coefficient on $D_L$ has the form   }
 \be\label{9}
\mu_{\wh f}(z) = \lambda_D^{-2}(z) \ov{\psi(z)}, \quad \psi = \tau(\vp).
\end{equation}

The Beltrami coefficients of such form are called now harmonic in view of their connection with the Kodaira-Spencer deformation theory of complex structures.

In this paper we shall use the notion of harmonicity in its original sense, i.e., only for solutions of the Laplace equation
$\Delta u = 0$, and regard the coefficients of type (9) as pseudo-harmonic.

\bigskip\bigskip
\centerline{\bf 2. GENERAL THEOREM AND ITS CONSEQUENCES}

\bigskip
Our goal is to describe the features of extremal Beltramic oefficients of non-Teichm\"{u}ller type maximizing the Grunsky norm.

We shall consider the univalent functions $f(z) \in S_Q(D^*)$ whose restrictions to the boundary quasicircle $L = \partial D^*$ have
substantional points, hence do not have the Teichm\"{u}ller extremal extensions across $L$. Denote the collection of
such $f$ by $S_Q^0(D^*)$. For $L = \mathbb S^1$, we use the notation $S_Q^0$.

Fix $p > 1$ and consider the subset $\mathcal M_p(D)$ of $\Belt(D)_1$, which consists of the Beltrami coefficients $\mu$
defining the maps $f^\mu \in S_Q^0(D)$ and satisfying:

$(i)$ \ $\mu \in L_\iy(D) \bigcap W^{1,p}(D)$,where $W^{1,p}(D)$ is the Sobolev space of functions $\mu$ in $D$ having the first
distributional derivatives which belong to $L_p(D)$;

$(ii)$ the value $\|\mu\|_\iy = \esssup_D |\mu(z)|$ is attained by approaching $z$ the boundary of $D$;

$(iii)$ there is a subarc $\gamma \subset \partial D$ depending on $\mu$ such that $\mu(z) \to 0$ as $z$ approaches $\gamma$
from inside $D$.

The boundary values $\mu(z_0)$ for $z_0 \in \partial D$ must be understand as $\lim\limits_{z \to z_0} \mu(z)$.
Note that the value  $\|\mu\|_\iy $ also can be attained by approaching the inner points of domain $D$ and that the arc $\gamma$ is
locally $C^\a$ smooth with $\a > 0$ depending on $\|\mu\|_\iy$, in accordance with the H\"{o}lder continuity of
quasiconformal automorphisms of $\hC$.

\bigskip
The main result of this paper is the following general theorem.

\bigskip\noindent
{\bf Theorem 1}. {\it For any $p > 2$, every Beltrami coefficient $\mu \in \mathcal M_p(D)$ is extremal in its equivalence class
$[\mu]$, and the corresponding quasiconformal automorphism $f^\mu$ of $\hC$ satisfies}
 \be\label{10}
k(f^\mu) = \vk_{D^*}(f^\mu) = \|\mu\|_\iy.
\end{equation}

\bigskip
In the case $p > 1$, there is a weakened version of this theorem, presented in Section 4.

Until now, there were no results in geometric function theory giving the explicit expression of the generalized Grunsky
norm for the general quasiconformal domains different from the disk.
Theorem 1 implies the first general result in this direction.

In view of the assumptions on the set $\mathcal M_p(D)$, every boundary function $f^\mu|\partial D$ must have at least one substantial point, at which $\|\mu\|_\iy$ is attained. This yields that the coefficient $\mu$ {\it is not uniquely extremal} in its equivalence class.

\bigskip
It follows from what was indicated above that the situation described by Theorem 1 does not appear in the case of sufficiently
high boundary regularity of univalent functions $f$ (and of $\partial D$) because, for example, for any  $C^{2+\a}$ smooth
$\mu$ the map $f^\mu$ is $C^{2+\a}$ on $\C$), and by Strebel's frame mapping criterion \cite{GL}, \cite{St2}) the equivalence
class of $f^\mu$ contains unique Teichm\"{u}ller coefficient $\mu_0 = k |\psi_0|/\psi_0$ with $\psi_0 \in A_1(D)$, which cannot
have the substantial boundary points.

The assumption on the set $\mathcal M_p$ are natural and cannot be replaced in terms of smoothness or non-smoothness of $\mu$
on the boundary points.

\bigskip
In the case, when the domain $D$ is the unit disk $\D$, Theorem 1 can be completed by the quantitative results on Fredholm
eigenvalues and reflection coefficients indicated in the previous section (see the relations (6), (7)).

\bigskip\noindent
{\bf Theorem 2}. {\it For any $p > 2$, every Beltrami coefficient $\mu \in \mathcal M_p(\D)$ is extremal in its equivalence class; in addition, the reflection coefficient and the Fredholm eigenvalue of the curve $L_1 = f^\mu(\mathbb S^1)$ are explicitly given by
 \be\label{11}
q_{L_1} = 1/\rho_{L_1} = \|\mu\|_\iy.
\end{equation}
}

Here also $\mu$ can be regular only in admissible bounds.
The assumption for $\mu$ to belong to the Sobolev space $W_1^p$ with $p > 2$ is essential for the proof.

In the case of pseudo-harmonic Beltrami coefficients, the assumptions of the above theorems can be weakened.

\bigskip\noindent
{\bf Theorem 3}. {\it The equalities (9) and (10) are valid for any admissible for univalence pseudo-harmonic Beltrami coefficient $\mu_\vp(z)$ of the form (9) or (8) with $\vp = S_f$ such that the maximum of the function $\ld_D^{-2}(z) \ov{S_f(z)} $ is attained
at some boundary point and $S_f$ is bounded on some subarc $\gamma$ of $\partial D^*$ .
}

\bigskip
For such coefficients, the required smoothness of $\mu$ is provided by the representations (8), (9), while vanishing on
$\gamma$ follows from (5). A special case of Theorem 3 has been obtained in \cite{Kr9}.

\bigskip\bigskip
\centerline{\bf 3. PROOF OF THEOREM 1}

\bigskip
To prove that the Beltrami coefficient $\mu$ satisfying the prescribed conditions $(i), (ii), (iii)$ is extremal in its equivalence class and provides the equalities (10), we establish that it must satisfy
 \be\label{12}
\|\mu\|_\iy = \sup_{{\|\psi\|}_{A_1^2(D)} = 1} \Big\vert \iint\limits_D \mu(z) \psi(z) dx dy \Big\vert.
\end{equation}
This equality means that $\|\mu\|_\iy$ is attained in $L_1(D)$ on the set of abelian quadratic differentials
intrinsically connected with the Grunsky coefficients.

By the Sobolev embedding theorem, the function $\mu(z)$ is extended to a continuous function on the closed domain $\ov D$.
The assumption $(ii)$ implies that there is a point $z_0 \in \partial D$ at which
 \be\label{13}
|\mu(z_0)| = \|\mu\|_\iy
\end{equation}
(in the general case,
$|\mu(z_0)| = \limsup\limits_{z \to \partial D} |\mu(z)|$).

We take the conformal map $z = \chi(\zeta)$ of the half-strip
$$
\Pi_{+} = \{\zeta = \xi + i \eta : \ \xi > 0, \ 0 < \eta < 1\}
$$
onto $D$ such that $\chi^{-1}(z_0) = \iy$ and the pre-images of the endpoints of the arc $\gamma$ are the points $0$ and $1$,
and pull back the coefficient $\mu$  to Beltrami coefficient
 \be\label{14}
\nu(\z) := \chi^* \mu(z) = (\mu \circ \chi)(\z) \ \ov{\chi^\prime(\z)}/\chi^\prime(\z)
\end{equation}
on $\Pi_{+}$. The horizontal lines
$$
l_\eta = \{\z = \xi + i \eta: \ 0 < \xi < \iy\}, \quad 0 < \eta < 1,
$$
are moving under this map into some analytic curves in $D$ with endpoints on $\partial D$. The infinite point $\z = \iy$ is substantial for the composite map $F^\nu = f^\mu \circ \chi$, and by assumptions on $\mu$ and from (14), we have
$$
\lim\limits_{\xi \to \iy} |\nu(\xi + i \eta)| = \|\nu\|_\iy = \|\mu\|_\iy.
$$
The assumptions on the coefficient $\mu(z)$ and the smoothness of conformal map $\chi$ imply that
$\nu \in L_\iy \bigcap W_p^2(\Pi_{+})$ and that the limit function
$$
\nu(\z_0) = \lim\limits_{\z \to \z_0 \in \partial \Pi_{+}} \nu(\z),
$$
is defined at all $\z_0 \in \partial \Pi_{+}$ different from the points $0, i, \iy$. Hence, $\nu(\z)$ is
bounded on the interval $[0, i]$ of the imaginary axes, and moreover,
 \be\label{15}
\nu(i\eta) = 0, \quad 0 \le \eta \le 1.
\end{equation}

Now we pick the sequence
$$
\omega_m(\z) = \frac{1}{m} e^{- \z/m}, \quad \z \in \Pi_{+} \ \ (m = 1, 2, ...);
$$
all these $\omega_m$ belong to $A_1^2(\Pi_{+})$ and $\omega_m(\zeta) \to 0$
uniformly on $\Pi_{+} \bigcap \{|\z| < M\}$ for any $M < \iy$. Also, $\| \omega_m \|_{A_1(\Pi_{+})} = 1$ (moreover,
$\bigl | \iint_{\Pi_{+}} \omega_m d\xi d\eta \bigr | = 1 - O(1/m)$), which shows that $\{\omega_m\}$
is a degenerating sequence for the affine horizontal stretching of $\Pi_{+}$.

We prove that this sequence is also degenerated for $\nu$, estimating the integrals
$$
I_m = \iint\limits_{\Pi_{+}} \nu(\z) \omega_m(\z) d \xi d \eta
$$
for large $m$. We have
  \be\label{16}
I_m = \int\limits_0^1 e^{-i \eta/m} d \eta \ \Bigl( \frac{1}{m} \int\limits_0^\infty \nu(\xi + i \eta) e^{- \xi/m}
d \xi\Bigr).
\end{equation}
The inner integral in (16) represents the values of the Laplace transform
$$
\mathcal L \nu = \int\limits_0^\iy \nu(t) e^{-st} dt
$$
of $\nu(\xi + i \eta)$ in the points $s = 1/m$, and the assumptions of the theorem imply the existence of this transform
also for the derivative
$\partial \nu(\xi + i \eta)/\partial \xi$, which is integrable over $\Pi_{+}$.

Integrating by parts and using (15), one obtains
 \be\label{17}
\int\limits_0^\iy \fc{\partial \nu(\xi + i \eta)}{\partial \xi} e^{- \xi/m} d \xi
= \fc{1}{m} \int\limits_0^\iy \nu(\xi + i \eta) e^{- \xi/m} d \xi - \nu(i \eta)
= \fc{1}{m} \int\limits_0^\iy \nu(\xi + i \eta) e^{- \xi/m} d \xi
\end{equation}
Since for any $s > 0$,
$$
\left | \fc{\partial \nu(\xi + i \eta)}{\partial \xi}\right | e^{- s \xi} <
\left | \fc{\partial \nu(\xi + i \eta)}{\partial \xi} \right |,
$$
the Lebesgue theorem on dominated convergence yields
$$
\lim\limits_{s \to 0} \int\limits_0^\iy \fc{\partial \nu(\xi + i \eta)}{\partial \xi} e^{-s \xi} d \xi
= \int\limits_0^\iy \fc{\partial \nu(\xi + i \eta)}{\partial \xi} d \xi.
$$
Together with (15) and (17), this implies
 \be\label{18}
\lim\limits_{m \to \iy} \fc{1}{m} \Bigl | \int\limits_0^\iy
\nu(\xi + i \eta) e^{- \xi/m} d \xi \Bigr | = |\nu(\iy) - \nu(i \eta)| = |\nu(\iy)|,
\end{equation}
where
$$
\nu(\iy) = \lim\limits_{\xi \to - \iy} \nu(\xi + i \eta).
$$
It follows that the integral (16) has the limit value
 \be\label{19}
\lim\limits_{m\to \iy} \iint\limits_{\Pi_{+}} \nu(\z) \omega_m(\z) d \xi d \eta
= \int\limits_0^1 d \eta \ \lim\limits_{m \to \infty} \fc{1}{m} \int\limits_0^\iy
\nu(\xi + i \eta) e^{- \xi/m} d \xi = \nu(\iy).
\end{equation}

Now we return to the initial domain $D$ by applying the inverse conformal map $\chi^{-1}(z): \ D \to \Pi_{+}$.
The corresponding sequence
$$
\psi_m = (\omega_m \circ \chi^{-1})(\chi^\prime)^{-2}, \quad m = 1,2, \dots,
$$
is a degenerating sequence for the initial Beltrami coefficient $\mu$ on $D$, and by (19),
 \be\label{20}
\lim\limits_{m\to \iy}
\Bigl | \iint\limits_D \mu(z) \psi_m(z) dx dy \Bigr |
= \lim\limits_{m\to \iy} \Bigl | \iint\limits_{\Pi_{+}} \nu(\zeta) \omega_m(\zeta) d \xi d \eta \Bigr | = |\nu(\iy)|.
\end{equation}
In view of the assumption (13), all terms in (20) are equal to $\|\nu\|_\iy = \|\mu\|_\iy$. This  proves the equality (12),
which implies that $\mu$ is extremal in its class and that the map $f^\mu$ obeys the relations (10).
The theorem follows.

\bigskip
Note that under the assumption $p > 2$, the conformal map $\chi$ of $\Pi_{+}$ onto $D$ is $C^{1,\a}$-smooth in the closed domain $\ov{\Pi_{+}}$, excluding the points $0, 1, \iy$, which yields that the restrictions $\nu_\eta$ of the Beltrami coefficient (14) to the lines $l_\eta$ belong to $W^{1,p}(l_\eta)$. Then the embedding theorem for $W^{n,p}$ functions giving the continuity of these restrictions implies the existence of the limits
$$
\nu(0) = \nu_\eta(0) = \lim_{\xi\to 0} \nu_\eta(\xi)\quad \text{for all} \ \ \eta,
$$
and $\nu_\eta(\iy) = \lim_{\xi\to \iy} \nu(\xi + i \eta)$.

This remark indicates the way in which the above proof can be extended to $p > 1$ (see remark {\bf 4.3}.

\bigskip\bigskip
\centerline{\bf 4. ADDITIONAL REMARKS TO THEOREMS 1 AND 2}

\bigskip\noindent
{\bf 4.1}. Theorem 1 implies that for every  Beltrami coefficient $\mu_0 \in \mathcal M_p(D)$ the disk
$$
\{t \mu_0/\|\mu_0\|_\iy: \ |t| < 1\}
$$
is geodesic in the ball $\Belt(D)_1$ simultaneously for the Teichm\"{u}ller, Kobayashi and Carath\'{e}odory metrics. It is pushed
down under the projection $\phi_\T: \ \mu \to S_{w^\mu} \in \B(D^*)$ onto a holomorphic disk in universal Teichm\"{u}ller space $\T$, on which all these  invariant distances also coincide,
though this coefficient does not be necessarily uniquely extremal in its equivalence class (the images of two different
$t^\prime \mu_0/\|\mu_0\|_\iy$ and $t^{\prime\prime} \mu_0/\|\mu_0\|_\iy$ in $\T$ can be the same; as well as, there are $\mu \in [\mu_0]$ with $\|\mu\|_\iy = \|\mu_0\|_\iy$).

\bigskip\noindent
{\bf 4.2}.
The equality (18) is equivalent to the Tauberian theorem for Laplace transform.
All known theorems of this type also require some smoothness of the original functions.

This nice connection with the Laplace transform implies simultaneously the extremality and equality of the Teichm\"{u}ller and
Grunsky norms norms for a broad set of Beltrami coefficients given by Theorems 1 and 2.

\bigskip\noindent
{\bf 4.3}. An extension of Theorem 1 to $p > 1$ is possible under the additional assumptions on $D$ and $\mu$,
for example, when the boundary of domain $D$ is $C^{1, \a}$-smooth ($\a > 0$) and in $(ii)$ the value
$\|\mu\|_\iy = \esssup_D |\mu(z)|$ is attained at some point $z_0 \in \partial D$ as $z \to z_0$ along any way in $D$,
then the above prove can be modified for $\mu \in L_\iy(D) \bigcap W^{1,p}(D)$ with $\ p > 1$ as follows.

Since now the conformal map $\chi$ of $\Pi_{+}$ onto $D$ also is $C^{1,\a}$-smooth in the closed domain $\ov{\Pi_{+}}$, excluding the points $0, 1, \iy$, we have that the restrictions $\nu_\eta$ of the Beltrami coefficient (14) to the lines $l_\eta$ belong to $W^{1,p}(l_\eta)$. Applying the embedding theorem in dimension $n = 1$, one obtains for $p > 1$ the existence of limits
$$
\nu_\eta(0) = \lim_{\xi\to 0} \nu_\eta(\xi) = 0 \quad \text{and} \ \ \nu_\eta(\iy) = \lim_{\xi\to \iy} \nu(\xi + i \eta).
$$
The indicated change of $(ii)$ yields that all values $\nu_\eta(\iy)$ must coincide and are equal to $\|\nu\|_\iy$.
Therefore, one can again apply the relations (18)-(20) and derive the assertions of Theorem 1.

\bigskip\noindent
{\bf 4.4}. The assertion of Theorem 1 can be extended to arbitrary quasiconformal domains, even with fractal boundaries
$L$ in the following form: one can take in the half-strip $\Pi_{+}$ the Beltrami coefficients $\mu$ obeying the Tauberian
theorem for the Laplace transform, i.e., for which the arguments of the proof of Theorem 1 are valid, and pull back
these $\mu$ to the interior of $L$.

\bigskip\noindent
{\bf 4.5}. The canonical quasiconformal extensions of univalent functions (with Teichm\"{u}ller or pseudo-harmonic of type (9)
coefficients $\mu$) are unique; the uniqueness intrinsically relates to their analyticity.

Theorem 1 involving Beltrami coefficients from Sobolev's space provides a possibility to distinguish some subsets in $\mathcal M_p$ admitting uniqueness.
For example, it holds for $\mu$ which are (even weak) solutions of the Dirichlet problem for uniformly elliptic differential
equations  of the second order $\mathcal L u = 0$ on $D$ with prescribed values $\mu$ on the boundary curve $\partial D$,
in particular, for harmonic $\mu(z)$ (with $\Delta \mu = 0$).
In addition, in this case the assumption (ii) is trivially fulfilled by the maximum principle.

In particular, all harmonic $\mu$ (with $\Delta \mu = 0$) vanishing on some boundary subarcs $\gamma$ are extremal, except
when the boundary function $f^\mu|\partial D$ satisfies the Strebel frame mapping condition (for example, it is $C^{2+\a}$ smooth).

\bigskip\noindent
{\bf 4.6}. The theory of extremal quasiconformal maps originated in \cite{Te} plays now a crucial role in quasiconformal analysis and in its deep applications to geometric function theory, Teichm\"{u}ller space theory and other  fields of mathematics and mathematical physics. The canonical extremal Beltrami differentials of Teichm\"{u}ller type
$\mu(z) d\ov z/dz$ with $\mu(z) = |\psi(z)|/\psi(z)$ generated by integrable holomorphic quadratic differentials $\psi(z) dz^2$
naturally arise in many problems. As was mentioned above, any such differential is unique in its equivalence class, and such maps $f^\mu$ are dense in $S_Q(D^*)$.  On unique extremality see also \cite{BLMM}.

The extremal quasiconformal extensions of univalent functions with equal Teichm\"{u}ller and Grunsky norms (hence, determined
by the squares of abelian differentials) have similar features.

The first example of extremal quasiconformal maps of non-Teichm\"{u}ller type was given by Strebel \cite{St1}. Recently,
the author found in \cite{Kr4} an important application of such coefficients  to geometric problems of Teichm\"{u}ller space theory. Other new types of not canonical extremal Beltrami coefficients are given in \cite{Kr9}.
The present paper continues this line, and Theorem 1 provides, in particular, that the structure of such coefficients can be rather pathological.

\bigskip\noindent
{\bf 4.7}. Since the quantities $k(f^\mu)$ and $\vk_{D^*}(f^\mu)$ depend continuously on $\mu$ and on the Schwarzians $S_{f^\mu}$
(respectively, in $L_\iy$ and $\B$ norms), all assertions of the above theorems are valid for the limit functions of sequences
$\mu_n \to \mu_0$ in the indicated norms.

\bigskip\bigskip
\centerline{\bf 5. GENERALIZATION}

\bigskip\noindent
{\bf 5.1. Improvement of Theorem 1}.
In the case, when the domain $D$ is the unit disk $\D$, we have the following strengthening of Theorems 1 and 2.
Its proof is much more complicated. It is based on an important theorem from \cite{Kr8} whose proof involves the deep results
on the Gaussian curvature, Grunsky inequalities and complex geometry of universal Teichm\"{u}ller space.

\bigskip\noindent
{\bf Theorem 4}. {\it Let the Beltrami coefficient  $\mu \in \Belt(\D)_1$ satisfy $\mu \in L_\iy(D) \bigcap W^{1,p}(\D)$
with some $p > 2$, and suppose that the value $\|\mu\|_\iy = \esssup_\D |\mu(z)|$ is attained by approaching $z$ the unit circle
$\mathbb S^1$ and on some subarc $\gamma$ of $\mathbb S^1$, we have
$$
\mu(z) \equiv q = \const \quad \text{with} \ \ |q| < \|\mu\|_\iy.
$$
Then $\mu$ is extremal in its class and the corresponding quasiconformal automorphism $f^\mu$ of $\hC$ admits the equalities}
\be\label{21}
k(f^\mu) = \vk(f^\mu) = q_{f^\mu(\mathbb S^1)} = 1/\rho_{f^\mu(\mathbb S^1)} = \|\mu\|_\iy.
\end{equation}

\bigskip\noindent
{\bf Proof}. Denote $D = f^\mu(\D), \ D^* = f^\mu(\D^*)$ and consider the maps $g^c$, which are conformal in $D^*$ and have in
$D$ a constant quasiconformal dilatation $c$. We regard such maps as the {\bf affine-like deformations of domain} $D$ and the collection of images $g^c(D)$ as the affine class of $D$. Each map $g^c$ has on $D$ the same Beltrami coefficient as the affine map
$$
\omega^c(w) = c_1 w + c_2 \ov w + c_3
$$
whose Beltrami coefficient with $\mu_{\omega^c} (w) = c_2/c_1 = c$ on $\C$ (so $g^c$ and $\omega^c$ differ on
$\omega \circ f^\mu(\D)$ by a conformal map).

Consider the map
$$
f_{c_0} (z) = g^{-q} \circ f^\mu(z),
$$
where $c_0 = - q$ coincides with the value of $\mu$ on the subarc $\gamma_0$ and is the Beltrami coefficient of the map
$g^{-q}(\omega)$ (inverse to affine deformation $g^q$).

By the chain rule for Beltrami coefficients, $f_q$ has the Beltrami coefficient
$$
\mu_{f_q}(z) = \fc{\mu(z) - q}{1 - q \ \ov{\mu(z)}} \ \fc{\partial_z f^\mu}{\ov{\partial_z f^\mu}}.
$$
It vanishes on the arc $f^\mu(\gamma) \subset \partial D$, and therefore the map $f_q$ satisfies the assumptions of
Theorems 1 and 2 on the disk $\D$, which imply for $f_q$ the corresponding equalities (10) and (11).

Now we apply the following theorem from \cite{Kr8}.

\bigskip\noindent
{\bf Theorem A}. {\it For any function $f \in \Sigma_Q$ with $\vk(f) = k(f) < 1$ mapping the disk $\D^*$ onto the complement
of a bounded quasidisk $D$ and any affine-like deformation $g^c$ of this domain (with $|c| < 1$), we have the equality}
$$
\vk(g^c \circ f) = k(g^c \circ f).
$$

\bigskip
Taking $f = f_q|\D^*$ and $c = q$, one obtains by Theorem A the equalities (21), completing the proof of Theorem 4.

\bigskip\noindent
{\bf 5.2. Geometric features}. If a Beltrami coefficient $\mu \in \mathcal M_p(D), \ p > 2$, is harmonic, then also all
$t \mu$ with $|t| < 1$ are harmonic. In view of their unique extremality, the image
$$
\D_h(\mu_0) = \{\phi_\T(t \mu): \ |t| < 1\}
$$
is a holomorphic disk (without singular points) in the universal Teichm\"{u}ller space $\T$. By Theorem 4 this disk
is geodesic in the Teichm\"{u}ller, Kobayashi and Carath\'{e}odory metrics on $\T$.
This improves the assertion of Remark {\bf 4.1}.

It was established in \cite{Kr5} that the Grunsky coefficients of univalent functions generate a Finsler structure
$G_\T(\vp, v)$  on the tangent bundle of the space $\T$, which is dominated by its canonical Finsler structure $F_\T(\vp, v)$
generating the Teichm\"{u}ller metric of this space.
The structure $G_\T(\vp, v)$ canonically generates the corresponding measurable infinitesimal Finsler metric $\lambda_\vk$,
and due to \cite{Kr5}, on any extremal Teichm\"{u}ller disk
$\D(\mu_0) = \{\phi_\T(t \mu_0): \ t \in\D\}$ and its isometric
images in $\T$, we have the equality
$$
\tanh^{-1}[\vk(f^{r\mu_0})] = \int\limits_0^r \ld_\vk(t) dt.
$$
The arguments applied in \cite{Kr5} remain in force also for harmonic geodesic disks $\D_h(\mu)$.

\bigskip\bigskip
\centerline{\bf 6. ILLUSTRATING EXAMPLES}

\bigskip
We illustrate the above theorems by the following examples.

\bigskip\noindent
{\bf Example 1: rectangles}. Even this case has been unsolved a long time.
Take a rectangle $P_4$ with vertices
$$
A = 0, \ \ B = a > 0, \ \ C = a + i b \ (b > 0), \ \ D = i b
$$
and define on
segments $0 \le x \le a$ and $0 \le y \le b$ two absolutely continuous functions $h_1(x)$ and $h_2(y)$ with the derivatives  $h_1^\prime(x), \ h_2^\prime(y) \in L_p, \ p >1$, and such that $h_1(0) = 0$ and both functions are not decreasing.
So, $h_1(x) h_2(y)$ vanishes on the left side of $P_4$. Assume also that $h_1(x) < 1, \ h_2(y) < 1$ and define
$$
\mu(z) = (1 + i) h_1(x) h_2(y), \quad z = x + iy.
$$
Any such $\mu$ satisfies the assumption of Theorem 1. The corresponding solution $f^\mu(z) = z + a_2 z^2 + \dots, \ |z| < 1$,
of the Beltrami equation $\partial_{\ov z} w = \mu(z) \partial_z w$ on $\C$ obeys the property (10).

Composing this function with conformal map of the disk $\D$ onto $P_4$ given by the Schwarz-Christoffel integral
$$
g(z) = \int\limits_0^z \fc{dt}{\sqrt{(t^2 - 1) (t - i) (t - \a)}},
$$
where the points $1, i, -1, \a$ are the preimages of the vertices of $P_4$ on $\mathbb S^1$, one obtains that
$f^\mu \circ \chi$ satisfies (10) and (11).

This strengthens the results of \cite{Kr3} on Fredholm eigenvalues of rectangles obtained by applying the Finsler geometry
of universal Teichm\"{u}ller space.

\bigskip\noindent
{\bf Example 2: ellipses}. Let $D$ is the ellipse with the foci at $-1, 1$ and semiaxes $a, b \ (a > b)$. An orthonormal basis in the Hilbert space $A_2(D)$ of the square integrable holomorphic functions on $D$ is formed by the polynomials
$$
P_n(z) = 2 \sqrt{\fc{n + 1}{\pi}} \ (r^{n+1} - r^{-n-1}) \ U_n(z),
$$
where $r = (a + b)^2$ and $U_n(z)$ are the Chebyshev polynomials of the second kind,
$$
U_n(z) = \fc{1}{\sqrt{1 - z^2}} \ \sin [(n + 1) \arccos z], \quad n = 0, 1, \dots
$$
(see \cite{Ne}).
Then for any Beltramu coefficient $\mu$ satisfying the conditions $(i), (ii), (iii)$, we have the equalities:
$$
\|\mu\|_\iy = \sup  \ \Big\{ \Bigl | \iint_D \mu(z) \sum\limits_0^\iy c_n P_n(z) dx dy \Bigr |: \
\left \| \sum\limits_0^\iy c_n P_n(z) \right\|_{A_1(D)} = 1 \Big\} = k(f^\mu) = \vk_{D^*}(f^\mu).
$$

\bigskip\noindent
{\bf Example 3: Analytic curvelinear polygons}. To simplify the formulas, we pass to quasiconformal automorphisms $f^\mu$ of $\hC$
conformal on the lower half-plane $H^* = \{z: \ \Im z < 0\}$ (instead of the disk). Their Beltrami coefficients $\mu$ are supported in the upper half-plane
$H = \{z: \Im z > 0\}$ and run over the ball
$$
\Belt(H)_1 = \{\mu \in L_\iy, \ \mu(z)|H^* = 0, \ \|\mu\|_\iy < 1\},
$$
and the Schwarzian derivatives $S_{f^\mu}$ belong to the space $\B(H^*)$ formed by holomorphic quadratic differentials
$\vp(z) dz^2$ on $H$ with norm $\|\vp\|_\B = \sup_{H^*} |z - \ov z|^2 |\vp (z)|$.

Pick unbounded convex rectilinear polygon $P_n$ with finite vertices $A_1, \dots, A_{n-1}$ and $A_n = \iy$. Denoting its exterior angles at $A_j$ by $\pi \a_j$ so that $\pi < \a_j < 2 \pi, \ j = 1, \dots, n - 1$, one obtains that the conformal map $f_n$
of the lower half-plane $H^* = \{z: \ \Im z < 0\}$ onto the complementary polygon $P_n^* = \hC \setminus \ov{P_n}$ is realized
by the Schwarz-Christoffel integral
$$
f_n(z) = d_1 \int\limits_0^z (\xi - a_1)^{\a_1 - 1} (\xi - a_2)^{\a_2 - 1} ... (\xi - a_{n-1})^{\a_{n-1} - 1} d \xi +
d_0,
$$
with $a_j = f_n^{-1}(A_j) \in \R$ and complex constants $d_0, d_1$; here $f_n^{-1}(\iy) = \iy$.
Its Schwarzian derivative equals
$$
S_{f_n}(z) = \mathbf b_{f_n}^\prime(z) - \fc{1}{2}  \mathbf b_{f_n}^2(z) =
\sum\limits_1^{n-1} \frac{C_j}{(z - a_j)^2} -
\sum\limits_{j,l=1}^{n-1} \frac{C_{jl}}{(z - a_j)(z - a_l)},
$$
where $\mathbf b_f = f^{\prime\prime}/f^\prime$,
$$
C_j = - (\a_j - 1) - (\a_j - 1)^2/2 < 0, \ \ C_{jl} = (\a_j - 1)(\alpha_l - 1) > 0.
$$
Denote by $r_0$ the positive root of the equation
$$
\frac{1}{2} \Bigl[\sum\limits_1^{n-1} (\alpha_j - 1)^2 +
\sum\limits_{j,l=1}^{n-1} (\alpha_j - 1)(\alpha_l - 1) \Bigr] r^2 -
\sum\limits_1^{n-1} (\alpha_j - 1) \ r - 2 = 0,
$$
and define
$$
S_{f_n,t} = t \mathbf b_{f_n}^\prime -  \mathbf b_{f_n}^2/2, \ t > 0.
$$

By Theorem 3 and Ahlfors-Weill,  every Schwarzian $r S_{f_n,r_0}$ with $0 <r < r_0$ generates a univalent function $w_r: H^* \to \C$ whose pseudo-harmonic Beltrami coefficient
$\nu_r(z) = - (r/2) y^2 S_{f_n,r_0}\ov z)$ in $H$ is extremal in its equivalence class, and
$$
k(w_r) = \vk(w_{r \circ \sigma}) = \fc{r}{2} \| S_{f_n,r_0}\|_{\B(H^*)},
$$
where $\sigma$ is the appropriate Moebius map of $\D^*$ onto $H^*$.

The point is that in view of extremality of pseudo-harmonic coefficients $\nu_r$ following from Theorem 3, the Schwarzians $S_{f^{\nu_r}}$ with $r > r_0$ close to $r_0$ cannot lie in the space $\T$ modelled by the Schwarzians; this relates to the well-known problem on starlikness of Teichm\"{u}ller spaces in Bers' embedding, cf. \cite{Kr5} and the references therein. 

The images $f^{\nu_r}(H)$ with $0 < r < r_0$ are curvilinear polygons with piecewise analytic boundaries (in particular,
spirals).

\bigskip\noindent
{\bf Example 4: Quasiconformal polygons with affine-like sides}.
We fix on the unit circle some points $a_1, a_2, \dots, a_n$ following counterclockwise and regard the disk $\D$ as a
polygon with vertices at these points. Take a function
$u(\theta) = c_1 e^{i \theta} + c_2 e^{- i \theta} + c_3$
for $\arg a_1 \le \theta \le \arg a_2$, with complex $c_1, \ c_2, \ c_3$ and $|u(\theta)| < 1$
and extend it to a function $\wh u(\theta)$ on $[-\pi, \pi]$ defining by Poisson integral a harmonic function $\mu(z)$ on the disk $\D$ with $|\mu(z)| < 1$, which belongs to $W^{1,p}, \ p > 2$, and such that its pull-back to the half-strip $\Pi_{+}$ is compatible
with the Tauberian theorem for the corresponding Laplace transform, giving the equality (20) with $\nu(\iy) = \|\mu\|_\iy$.
We continue this $\mu$ by zero to $\D^*$.

By Theorem 4, the coefficient $\mu$ is extremal in its class, and the corresponding quasiconformal homeomorphism $f^\mu$ of $\hC$ obeys the relations (21).

\medskip
{\small \em{ \leftline{Department of Mathematics, Bar-Ilan University, 5290002 Ramat-Gan, Israel}
\leftline{and Department of Mathematics, University of Virginia, Charlottesville, VA 22904-4137, USA}
}}

\end{document}